\documentclass[5p,times]{elsarticle}

\usepackage{lineno,hyperref}
\usepackage{color}
\usepackage{verbatim}
\usepackage{bm}
\usepackage{amsmath}
\usepackage{graphicx}
\usepackage{soul}
\usepackage{multicol}
\usepackage{eqparbox}
\usepackage{algorithmic}

\modulolinenumbers[1]
\usepackage[ruled]{algorithm2e}
\usepackage{mathtools}

\newcommand{\Ref}[1]{\mbox{(\ref{#1})}}

\def\bq{{\mathbf q}}
\def\bv{{\mathbf v}}
\def\bd{{\mathbf d}}
\def\br{{\mathbf r}}
\def\bb{{\mathbf b}}
\def\bx{{\mathbf x}}
\def\bw{{\mathbf w}}
\def\b0{{\mathbf 0}}
\def\bu{{\mathbf u}}
\renewcommand{\div}{\operatorname{div}}
\newcommand{\R}{{\mathcal R}}
\newcommand{\M}{{\mathcal M}}

\journal{Parallel Computing}

\begin{document}

\begin{frontmatter}

\title{AMG based on compatible weighted matching\\ for GPUs}
\tnotetext[AMG for GPUs]{This work is partially supported by the EC under the Horizon 2020 Project
{\em Energy oriented Centre of Excellence for computing applications - EoCoE}, Project ID: 676629.}

\author[address1]{Massimo Bernaschi}
\author[address2]{Pasqua D'Ambra}
\author[address1,address3]{Dario Pasquini}
\cortext[mycorrespondingauthor]{Dario Pasquini}

\address[address1]{Institute for Applied Computing (IAC) - CNR, dei Taurini 19, 00185 - Rome, I.}
\address[address2]{Institute for Applied Computing (IAC) - CNR, Naples branch, Via P. Castellino, 111, 80131, Naples, I.}
\address[address3]{Department of Computer Science, ``Sapienza'' University -, Via Salaria, 113, 00198, Rome, I.}

\begin{abstract}
We describe main issues and design principles of an efficient implementation, tailored to recent generations of Nvidia Graphics Processing Units (GPUs), of an Algebraic MultiGrid (AMG) preconditioner previously proposed by one of the authors and already available in the open-source package {\em BootCMatch: Bootstrap algebraic multigrid based on Compatible weighted Matching} for standard CPU.  The AMG method relies on a new approach for coarsening sparse symmetric positive definite (s.p.d.) matrices, named {\em coarsening based on compatible weighted matching}. It exploits maximum weight matching in the adjacency graph of the sparse matrix, driven by the principle of compatible relaxation, providing a suitable aggregation of unknowns which goes beyond the limits of the usual heuristics applied in the current methods. We adopt an approximate solution of the maximum weight matching problem, based on a recently proposed parallel algorithm, referred as the \emph{Suitor algorithm}, and show that it allow us to obtain good quality coarse matrices for our AMG on GPUs. We exploit inherent parallelism of modern GPUs in all the kernels involving sparse matrix computations both for the setup of the preconditioner and for its application in a Krylov solver, 
outperforming preconditioners available in
Nvidia AmgX library.  We report results about a large set of linear systems arising from discretization of scalar and vector partial differential equations (PDEs).
\end{abstract}


\end{frontmatter}


\section{Introduction}
\label{intro}

We are concerned with efficient solution, on recent generations of GPU accelerators, of systems of linear equations:
\begin{equation}
A \bx=\bb,
\label{sys1}
\end{equation}
where $A \in \R^{n \times n}$ is a symmetric positive definite (s.p.d.), large and sparse matrix. More specifically, we focus on main issues and design principles driving a parallel implementation of all the functionalities of the package {\em BootCMatch: Bootstrap algebraic multigrid based on Compatible weighted Matching}~\cite{bootcmatch}, for preconditioning and solving system~\Ref{sys1} by an Algebraic MultiGrid (AMG) method based on aggregation.

AMG methods are a popular choice for dealing with a system like ~\Ref{sys1}, when it results from the discretization of partial differential equations (PDEs) on complex geometries and unstructured grids or when no information about its origins are available. Main distinguish feature of the above methods, with respect to their geometric counterpart, is the chance of defining an automatic setup of the hierarchy of coarse-level variables and matrices by relying only on the (fine) coefficient matrix. Many variants of AMG methods~\cite{RS1987,VMB1996,HY2000,N2010} and related parallel software libraries~\cite{FJY2006,ml-guide,agmg-guide} have been proposed in the literature. They differ in the way in which coarse-level variables are selected and in the setting of coarse-to-fine transfer operators. Despite of differences, current AMG methods show good numerical scalability, meaning that the number of iterations stays almost constant while the system size scales up, when they are applied to classes of sparse matrices corresponding to discretizations of 2nd order scalar elliptic PDEs.

In~\cite{DV2013,DFV2018} the authors propose a new AMG method, and the corresponding sequential software, which relies on a new setup procedure to generate coarse-level
variables
aimed at obtaining an
AMG preconditioner showing good numerical scalability for more general s.p.d. linear systems.

Here we present a parallel version of BootCMatch, which efficiently exploits the fine-grained parallelism and the memory organization of modern GPU accelerators, with the final aim to move a step towards AMG for future exascale computations.
In the last 10 years a growing number of systems embed GPU accelerators to exploit their outstanding performance (see~\cite{top500}). However, these new platforms may require to rethink and redesign algorithms, data structures and software development paradigms for taking full advantage from their usage. In particular, it is not unusual for algorithms that on traditional computing platforms are considered inefficient due to slow convergence, to become very much competitive on GPUs since additional computations are well tolerated and convenient with respect to using complex memory access patterns. For example, it is well known that highly parallel smoothers, such as versions of weighted Jacobi, outperform, in terms of execution times, more effective but intrinsicallly sequential smoothers as Gauss-Seidel relaxation for preconditioning and solving sparse linear systems on GPUs (e.g., see~\cite{DF2016,GEZ2014,BCHZ2014}). One of the main objectives of our work, as described
  in Section~\ref{gpu}, has been to implement highly tuned kernels that access GPU global memory according to best practices of CUDA programming and to use the available
  computing resources ({\em i.e.,} CUDA cores) in a cost-effective way by introducing the concept of {\em miniwarp}, for all the kernels implemented in BootCMatch.
The rest of the paper is organized as follows: Section~\ref{backg}
introduces the AMG methods and in particular the
variant that relies on the solution of a weighted graph matching
problem for the generation of the coarse-level variables. Section \ref{gpu} describes the issues related to
a parallel implementation of the AMG method based on compatible weighted matching. Section \ref{relwork} provides
a brief description of related works. Section \ref{results} presents the results obtained on a large set
of test cases. Finally, Section \ref{concl} concludes the work presenting future lines of activity.

\section{Background}
\label{backg}

\subsection{Algebraic Multigrid Methods}
\label{amg}

Multigrid methods are linear complexity methods for solving system~\Ref{sys1}. They are built on a relaxation method (the \emph{smoother}), such as a Ri\-chardson-type method, which efficiently damps high-frequency errors, although it is not able to reduce low-frequency errors. However, moving the problem to a coarser grid, what were previously low-frequency errors become high-frequency errors and can be damped by a new application of relaxation. The above procedure, whose setup requires coarser grids and transfer operators for moving among the grids, can be recursively applied obtaining methods with a computational cost which depends only linearly on the problem size~\cite{BHM2000,V2008}. While geometric multigrid methods rely on a pre-defined hierarchy of grids and on transfer operators depending on the geometry of the problem, AMG methods use only the information available in the system matrix. In the following, we describe the main components for setup and application of an AMG method; for an exhaustive introduction to AMG we refer the reader to~\cite{S2001}.

Let the set of row indices of $A$ be the fine index space, i.e., $\Omega = \{1, 2, \ldots, n\}$.
Any AMG generates a hierarchy of $nl$ index spaces and a corresponding hierarchy of matrices,
\[ \Omega^1 \equiv \Omega \supset \Omega^2 \supset \ldots \supset \Omega^{nl},
\quad A^1 \equiv A, A^2, \ldots, A^{nl}, \] by a suitable \emph{coarsening algorithm} using the information contained in $A$.
A vector space $\R^{n_{k}}$ is associated with $\Omega^k$,
where $n_k$ is the size of $\Omega^k$.
For all $k < nl$, a prolongation operator is built $P^k \in \R^{n_k \times n_{k+1}}$ and the matrix $A^{k+1}=(P^k)^TA^kP^k$ is computed according to the Petrov-Galerkin approach. A smoother operator $M^k$ is also defined, representing the iteration matrix of a relaxation method. All the above components are built in the so-called \emph{setup phase}.

The components produced in the setup phase may be combined in several ways
to obtain different types of \emph{multigrid cycles};
this is done in the \emph{application} or \emph{solve phase}. An example of such a combination, known as
symmetric V-cycle, is given in Algorithm~\ref{Vcycle_alg}. In that case, a single iteration
of the same smoother is used before and after the recursive call to the V-cycle (i.e.,
in the pre-smoothing and post-smoothing phases). However, more robust, although more expensive, choices can be
performed, such as W-cycle~\cite{BHM2000}, and recursive Krylov-based cycle (K-cycle)~\cite{NV2008}. At the coarsest level, i.e., for $k=nl$, a direct solver is usually employed. Actually, especially in parallel implementations of the algorithm, an iterative solver of the coarsest system is also applied in order to reduce data dependencies among parallel processors.
{\small
\begin{algorithm}[h]
\SetAlgoLined
V-cycle($k,A^k,\bb^k,\bx^k$)\\
\eIf{ $k \ne nl$}
{
$\bx^k = \bx^k + (M^k)^{-1} \left(\bb^k - A^k \bx^k\right)$\;
$\bb^{k+1} = (P^{k+1})^T\left(\bb^k - A^k \bx^k\right)$\;
$\bx^{k+1} =$ V-cycle$\left(k+1,A^{k+1},\bb^{k+1},\b0\right)$\;
$\bx^k = \bx^k + P^{k+1} \bx^{k+1}$\;
$\bx^k = \bx^k + (M^k)^{-T} \left(\bb^k - A^k \bx^k\right)$\;
}
{
$\bx^k = \left(A^k\right)^{-1} \bb^k$\;
}
\Return $x^k$
\caption{\emph{V-cycle}}
\label{Vcycle_alg}
\end{algorithm}
}

The choice of the coarse index spaces and of the prolongation operators are strictly related each other and affects the convergence properties of Algorithm~\ref{Vcycle_alg}, which, in turn, strongly depend on the ability of the coarse vector spaces to well represent the errors unaffected by relaxation (\emph{algebraically smooth vectors}) and of the prolongators to well interpolate them back to the fine space. Recent theoretical developments provide general approaches to the construction of coarse spaces for AMG having optimal convergence, i.e., a convergence independent of the problem size, in the case of general linear systems (see~\cite{XZ2017} and the references herein). However, despite these theoretical developments, almost all currently available AMG methods and software rely on heuristics to drive the coarsening process among variables; for example the strength of connection heuristics is derived from a characterization of the algebraically smooth vectors that is theoretically well understood only for M-matrices. The above heuristics is generally used both in the \emph{classical coarsening} and in an alternative approach, named \emph{coarsening by aggregation}~\cite{BHM2000}. The classical coarsening separates the original index set into either coarse indices (C-indices), which form the coarse level, and fine indices (F-indices), whose unknowns will be interpolated by the C-indices values, while aggregation-based coarsening uses aggregates of fine indices to form the coarse indices. In~\cite{DV2013,DFV2018} a new coarsening algorithm, which does not require a priori characterization of smooth vectors, has been proposed. It relies on the so-called {\em compatible relaxation} principle
introduced in~\cite{B2000}, which indicates a way to measure the quality of a coarse-level space, and exploits a maximum weight matching in the graph defined by the system matrix to find out an automatic aggregation-based coarsening for general s.p.d. matrices. In the following we describe the main features of the above aggregation algorithm and refer the reader to the original papers for details on the rationale and numerical principles at the base of its use for efficient coarsening.

\begin{figure*}[h]
\centering
\includegraphics[width=0.8\textwidth]{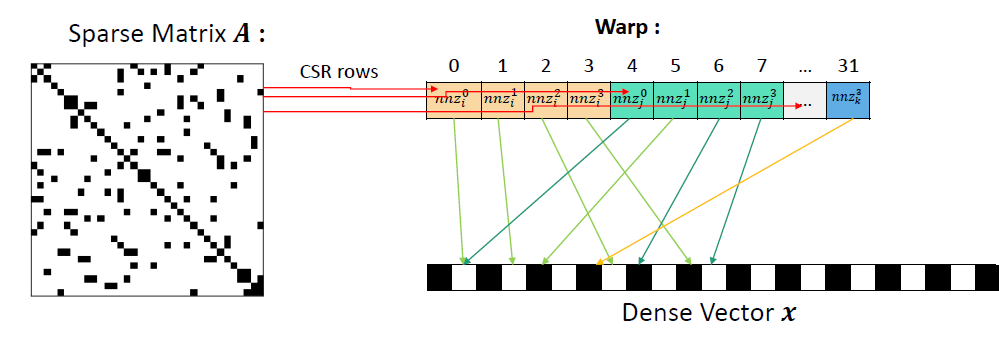}
\caption{Each {\em miniwarp} is in charge of a row of the matrix stored in CSR format.\label{fig-miniwarp}}
\end{figure*}

\subsection{Aggregation algorithm based on weighted graph matching}
\label{AMGmatching}

Let $G=(V,E,C)$ be the weighted undirected adjacency graph of the matrix $A$ in~\Ref{sys1}, where the vertex set $V$ consists of the row/column indices of $A$, the edge set $E$ corresponds to the couples of indices $(i,j)$ of the nonzero entries in $A$, and $C=(c_{ij})_{(i,j) \in E}$ is a matrix of positive edge weights.
A \emph{matching} in $G$ is a subset of edges $\M \subseteq E$ such
that no two edges share a vertex.
A {\em maximum weight matching} in the graph $G$ is defined as the $\arg \max_{\M} \sum_{ (i,j) \in \M} c_{ij}$.

In~\cite{DV2013,DFV2018} a maximum weight matching has been exploited to form aggregates of index pairs for good-quality coarsening in AMG methods. Main element driving the aggregation scheme is the definition of a suitable matrix $C(A,\bw)$ of edge weights for the adjacency graph of the original system matrix which is function of $A$ and of a vector $\bw \in \R^n$. More specifically, it consists of positive values arising from a linear-complexity computation involving the entries of matrix $A$ and of a good sample $\bw$ of algebraically smooth vectors for the system at hand. In principle, the vector can be arbitrary; a constant vector being a possible choice, however, in the BootCMatch framework it is a user-defined parameter.
The basic pairwise aggregation algorithm is described in Algorithm~\ref{alg-aggr}.
{\small
\begin{algorithm}[h]
\SetAlgoLined
\KwData{$G=(V,E,C(A,\bw))$, weighted adjacency graph of $A$}
\KwResult{$n_p$, $n_s$, $n_c$ and sets of aggregates $e_1, \ldots e_{n_c}$}
\begin{itemize}
\item Compute $\M$ \textbf{maximum weight matching} for $G$.
\item Initialize:
$n_c=0$, $n_p=0$, $n_s=0$\;
$U=[1, \ldots, n]$\;
\item
\While{$U \neq \emptyset$}
{
  Pick $i \in U$\;
\eIf{$\exists j \in U\setminus\{i\}$ such that $(i,j) \in \M$}
{
$n_p=n_p+1$\;
$n_c=n_c+1$\;
$e_{n_c}=\{i,j\}$\;
$U = U\setminus\{i,j\}$\;
}
{
$n_s=n_s+1$\;
$n_c=n_c+1$\;
$e_{n_c}=\{i \}$\;
$U = U\setminus\{i\}$\;
}
}
\end{itemize}
\caption{\emph{Pairwise aggregation based on maximum weight matching}}
\label{alg-aggr}
\end{algorithm}
}
Once the aggregates are formed, a piecewise constant interpolation operator, obtained by the orthonormal projection of the smooth vector $\bw$ on
the aggregates, is defined for the construction of a multigrid hierarchy in a recursive scheme. We note that in the recursive application of Algorithm~\ref{alg-aggr}, the input weighted graph $G=(V,E,C(A,\bw))$ corresponds to the adjacency graph of the computed coarse matrix whose weights are obtained by involving the restriction of the original vector $\bw$ on the coarse space.
Finally, we observe that multiple steps of the basic pairwise aggregation can be combined for obtaining larger aggregates in more aggressive coarsening.

Accurate solutions for the computation of maximum weight matching in a graph are based on the Hungarian algorithm to search optimal augmenting paths in the matrix between unmatched vertices~\cite{DK2001}. That algorithm has a super-linear worst-case complexity and it is intrinsically sequential, therefore it represents the main issue in the search for an efficient parallel computation of a maximum weight matching. However, approximate solutions featuring near-linear complexity have been shown to represent a viable approach to obtain a more efficient matching for good-quality coarsening and some of them have been included in the BootCMatch software framework.

\section{Parallel Algorithms for GPUs}
\label{gpu}

In the following we describe the design principles applied in our implementation, specifically tailored for recent generations of Nvidia GPUs using the CUDA framework, of the main kernels involved in setup and application of the compatible weigh\-ted matching AMG procedure as preconditioner in a precon\-di\-tio\-ned Conjugate Gradient (CG) method. As in the original BootCMatch code and in the Nvidia AmgX library, we chose to employ a CSR (Compressed Sparse Row) storage format for the sparse matrices.

\subsection{Setup of the preconditioner}
\label{build}

The two main issues in the development of a GPU version of the setup phase are: the computation of maximum weight matching in weighted graphs and the computation of coarse matrices.

The original CPU version offers the choice among a set of optimal and approximate maximum weight matching algorithms that, however, are either
inherently sequential or unsuitable to a good GPU implementation.
Therefore, we employed a different algorithm, named {\em Suitor}, which is a near-optimal matching algorithm recently proposed for GPU~\cite{HFVTP2012} and available in source form. We note that the CUDA kernels in Suitor original implementation make use of the {\em shuffle} instructions, a feature available starting from the {\em Kepler} architecture that offers a way to directly share data among threads belonging to the same warp (a group of 32 threads).
The original version of warp-level primitives depended on implicit warp-synchronous behaviour that, however, is not longer guaranteed starting on CUDA 9.0. Therefore, we adapted the Suitor algorithm to the new thread-scheduling policy supported by CUDA 9.0, although CUDA supports, by using suitable compiler options, a behaviour of the warp compatible with the legacy environment. We tested both options ({\em i)} updating Suitor so that it uses the new shuffle primitives (with explicit synchronization) and {\em ii)} using compiler options to use the legacy warp (synchronous) behaviour, but we did not find significative differences in the execution time of our code. We observe that, although the Suitor algorithm computes an approximate matching, it makes possible to obtain good quality coarsening in our AMG, as shown in the results presented in Section~\ref{results}.

After optimizing parallel execution of Algorithm~\ref{alg-aggr}, as also remarked in related work (see Section~\ref{relwork}), the most time-con\-su\-ming com\-pu\-ta\-tion remains the triple-matrix product involved in the Petrov-Galerkin approach for computation of coarse matrices. In the beginning, we resorted to the standard kernel available in the {\em cusparse}~\cite{CUSPARSE} library ({\em cusparseDcsrmm}).
However, we found that its performance was far from being optimal and we changed our code to use {\em Nsparse}, a recent implementation of sparse matrix-matrix product available in open source format~\cite{NNM2017}. Nsparse, as the implementation of Suitor, relies on the legacy shuffle primitives, nevertheless it provides a clear advantage with respect to the general-purpose primitives available in {\em cusparse}.

All other matrix operations of the setup phase, that are the computation of the transpose of the prolongator and the restriction of the current-level smooth vector, are implemented by assigning one instance of what we call {\em miniwarp} to each row of the matrix. A miniwarp is
a subset of a full warp (a set of $32$ threads) that consists of $\{2, 4, 8, 16\}$ threads. The choice of the size of the miniwarp depends
on the average number of nonzero entries per row of the sparse matrix. The main advantage of the miniwarp is that, for matrices with
few nonzero entries per row, the number of idle threads decreases. With the full warp, if a row has, on average, only $k<32$ nonzero entries,
there are, always on average, $32-k$ threads that remain idle. The mini warp reduces the difference significantly by using a size that is much closer to the average number of nonzero entries per row. Fig.~\ref{fig-miniwarp} shows how the miniwarp works when applied to a matrix in CSR format.

\subsection{Application of the preconditioner}
\label{solve}

In the preconditioner application, within the solve phase, we focused on two main kernels: the application of the smoothing in the multigrid cycle and the implementation of an optimized version of the CG method.
We note that in our code we implemented a flexible version of the CG method~\cite{NV2008}, needed in the case of variable preconditioners, such as Krylov-based AMG cycles (K-cycle).

For smoothing, we chose a version of Jacobi relaxation already used in~\cite{BCHZ2014} for AMG in a GPU setting. It is the so-called $\ell^1-$Jacobi smoother, which is a paramater-free version of the highly parallel Jacobi method always convergent for s.p.d. matrices and having good smoothing properties for strictly diagonally dominant matrices. The same kernel is used also for the solution of the linear systems at the coarsest level. Our implementation of this kernel relies again on miniwarps, each miniwarp is in charge of a row and the selection of the miniwarp size follows the same criterion above defined. The most expensive computation both for $\ell^1-$Jacobi smoother and for the preconditioned CG is the product between a sparse-matrix and a dense-vector that we indicate with {\em SpMV}, therefore we
focused on the tuning of that kernel for our aims. The sparse matrix involved in a {\em SpMV} can be:
\begin{itemize}
\item a coarse matrix; in this case, no assumption can be done on the number of nonzero entries per row;
\item a prolongator; in our aggregation scheme, also known as \emph{plain or unsmoothed aggregation}, the matrix has a single nonzero entry per row;
\item a transposed prolongator; the matrix has a number of nonzero entries per row that is, at most, equal to the size of the aggregates.
\end{itemize}
For the first case, depending on the sparsity degree of the matrix, the product is implemented  by using either a custom kernel that relies
on the concept of miniwarp or by the general-purpose {\em cusparse} primitive for the {\em SpMV}.
More precisely, the {\em cusparse} kernel is used when the number of non-zero entries per row is, on average, at least equal to the number of threads in a full warp (i.e., $32$ threads).
Indeed, we observed that, when the input matrix tends to be more dense, the {\em cusparse} primitive performs better.
The same technique is used for the transposed prolongator matrix, but in that case, the miniwarp product perfectly fits with the prolongator sparsity pattern.

For the prolongator matrix, it is possible to execute the {\em SpMV} more efficiently taking into account that the matrix has a single nonzero entry per row. In that case we used, for any row, a single thread.
The miniwarp approach employed in the {\em SpMV}, compared to the first {\em cuSPARSE} based implementation, provides, on average, a $1.4 \times$ speedup of the total solving time. That speedup increases up to $2.5 \times$ when the building phase produces a hierarchy composed by matrices with a low variance of $nnz$ per row.

For an efficient implementation of the preconditioned CG method, besides optimization of the SpMV computations, we also focused on reducing the number of GPU global memory access operations by employing a version of CG (see Algorithm~\ref{PCG}), originally proposed for a distributed implementation in~\cite{NN2015}. It computes a sequence of three scalar products within the main loop, by using a single kernel.
This approach allows us to reduce of a factor three the number of memory access operations, since the values of the vector $\bw$ that is involved in all the three products can be maintained in the registers. However, the reordering proposed in~\cite{NN2015} requires an additional {\em AXPY} computation(the update of a vector $Y$  as $Y=\alpha X + Y$). We grouped the {\em AXPY} computations in two pairs that are executed in a single kernel. The results of the first {\em AXPY} are maintained in the GPU registers and used for the second {\em AXPY} operation so reducing, again, the number of GPU global memory access operations. The optimized preconditioned CG implemented in our code is described in Algorithm~\ref{PCG}.
\begin{algorithm}
\small{
\label{PCG}
\caption{Preconditioned Flexible Conjugate Gradient}
\begin{algorithmic}[1]
\STATE Given $\bu_0$ and set $\br_0 = \bb - A\bu_0$
 \STATE $\bw_0 = \bd_0 = \mathcal{B}(\br_0)$
\STATE $\bv_0 = \bq_0 = A \bw_0$
\STATE $\alpha_0 = \bw_0^T \br_0$
\STATE $\beta_0 = \rho_0 = \bw_0^T \bv_0$
\STATE
	\STATE $\bu_{1} = \bu_0 + \alpha_0 / \rho_0 \bd_0$ 	
	\STATE $\br_{1} = \br_0 + \alpha_0 / \rho_0 \bq_0$
\STATE

\FOR{$i = 1,\dots$}
	\STATE $\bw_{i} = \mathcal{B}(\br_{i})$
	\STATE $\bv_{i} = A \bw_{i}$
	\STATE
	\STATE $\alpha_i = \bw_{i}^T \br_{i}$
	\STATE $\beta_i = \bw_{i}^T \bv_{i}$
	\STATE $\gamma_i = \bw_{i}^T \bq_{i-1}$
	\STATE $\rho_i = \beta_i - \gamma_i^2 / \rho_{i-1}$
	\STATE
	\STATE $\bd_{i} = \bw_{i} - \gamma_i / \rho_i \bd_{i-1}$
	\STATE $\bu_{i+1} = \bu_i + \alpha_i / \rho_i \bd_{i}$
	\STATE
	\STATE $\bq_{i} = \bv_{i} - \gamma_i / \rho_i \bq_{i-1}$
	\STATE $\br_{i+1} = \br_i + \alpha_i / \rho_i \bq_{i}$
	\STATE
\ENDFOR
\end{algorithmic}
}
\end{algorithm}

\section{Related Works}
\label{relwork}

Preconditioning and solving ever more large, sparse linear systems is a key kernel in computational science and the need to exploit the potential of GPUs for parallel preconditioners in iterative linear solvers is widely recognized~\cite{LS2013}. Due to their flexibility and potential scalability, many efforts were in particular devoted to parallel versions of AMG preconditioners specifically tailored to use single and multiple GPUs.

Some works~\cite{WRW2012,KF2012,ELB2012,LYC2015} focused on benchmarks of well-known AMG algorithms, such as AMG based on classical C/F coarsening~\cite{RS1987,HY2000} and aggregation-based AMG~\cite{VMB1996,N2010}, by using GPUs to accelerate the application of the preconditioner at each iteration of Krylov methods. They rely on efficient implementations of the SpMV kernel and emphasize that the setup of an AMG is a bottleneck in parallel AMG methods due to the sequential nature of the coarsening processes. On the other hand, focusing on accelerating application phase of AMG is justified by the need to repeat the above application iteratively in a Krylov process. Furthermore, it is frequent the need of solving many linear systems with the same matrix but different right-hand sides, e.g., in time-dependent or in Newton-type methods, therefore the setup cost can be amortized by multiple application phases.

Early work devoted to obtain a GPU implementation of both AMG setup and application phase on a single GPU is presented in~\cite{BDO2012}.
The authors describe main issues and their choices in implementing a version of the smoothed aggregation-based AMG proposed in~\cite{VMB1996}.
They rely on a fine-grained parallel implementation of a generalized maximal independent set algorithm for producing aggregates with similar properties and focus on efficient kernels for the Galerkin triple-matrix multiplication which represents the main roadblock on the way to obtain efficient AMG setup.

In~\cite{BCHZ2014} the authors present a GPU implementation of an unsmoothed ag\-gre\-ga\-tion-based AMG, where the focus
is both to implement an efficient parallel algorithm for computation of maximal independent set of variables, specifically tuned for standard isotropic graph Laplacian arising in 2nd order elliptic PDEs, and to simplify the Galerkin triple-matrix multiplication. Indeed, when standard unsmoothed aggregation is employed, the prolongation operator is a binary matrix and the Galerkin multiplication is reduced to summations of entries in the matrix at the finer level, which can be efficiently implemented in CUDA. They also emphasize that more sophisticated cycles than the standard V-cycle, such as K-cycle, should be employed in the case of unsmoothed-type aggregation schemes in order to preserve optimal convergence of the multilevel AMG method. Furthermore, they propose to use the $\ell^1-$Jacobi method both as smoother and as coarsest solver.

An efficient implementation for GPU of an unsmoothed\\ aggregation-based AMG is also discussed in~\cite{GEZ2014} and it is at the base of the GAMPACK commercial code~\cite{gampack} running both on single and multiple GPUs. Also in that work, the authors rely on a parallel algorithm for maximal independent set of coarse variables, representing aggregates of strictly connected fine variables, and propose to use a hybrid cycle to accelerate convergence of the application phase. They use K-cycle at the first 2 levels of the AMG hierarchy, whereas V-cycle is employed at the successive levels, in order to obtain a tradeoff between parallel efficiency and optimal convergence.

A description of the algorithms included in the publicly available Nvidia AmgX library~\cite{amgxsw}, running on single and multiple-GPUs, is in~\cite{AMGX2015}. AmgX implements both classical and unsmoothed aggregation-based AMG methods, with different choices for coarsening and prolongation operators. The parallel implementation of classical AMG is largely based on the methods implemented in the last available version of the hypre library~\cite{FJY2006}.
The aggregation algorithm of AmgX is based on a pairwise scheme similar to that proposed in~\cite{N2010}, coupling strongly-connected variables, which relies on a parallel graph matching techniques for efficient coarsening on single GPU. The library makes available a variety of cycles, such as V and W, and smoothers and coarsest solvers, including weighted-Jacobi, $\ell^1-$Jacobi, block-Jacobi, Gauss-Seidel, and an incomplete-LU (ILU) factorization. AmgX is the state of the art of AMG preconditioners for GPUs and in this paper we consider its single-node version for our performance comparisons.

\section{Numerical Experiments}
\label{results}

Hereafter, we discuss results obtained by using our GPU version of BootCMatch, named \emph{BootCMatchG}, for the solution of linear systems arising from scalar and vector PDE problems, as explained in the following.
For the first and second test cases, we consider linear systems of increasing size, in order to analyze scalability of the parallel preconditioners, whereas for the third case, we fix dimension and increase anisotropy, to the purpose of analyzing also the robustness of the methods.
\begin{description}
\item[ANI] These test cases derive from the following anisotropic 2D PDE on the unit square, with homogeneous Dirichlet boundary conditions:
\[
- \div(K\; \nabla u)=f,
\]
where $K$ is the coefficient matrix
\[
K = \left [
\begin{array}{ll}
a & c\\
c & b
\end{array}
\right ], \quad
\text{ with } \quad \left \{
\begin{array}{l}
a= \epsilon + \cos^2(\theta)\\
b= \epsilon + \sin^2(\theta)\\
c=  \cos(\theta)\sin(\theta)
\end{array}
\right .
\]
The parameter $0 < \epsilon \leq 1$ defines the strength of anisotropy in the problem,  whereas the parameter $\theta$ specifies the direction
of anisotropy. In the following we discuss results related to test cases with $\epsilon=0.001$ and $\theta= 0$, $\pi/8$, which we refer to as $ANI1$ and $ANI2$, respectively.
The problem was discretized using the Matlab PDE toolbox, with linear finite elements on (unstructured) triangular meshes of three different sizes ($168577$, $673025$, $2689537$), obtained by uniform refinement.
\item[LE] A second set of test cases comes from the discretization of the following Lam\'e
equations for linear elasticity:
\[
\mu \Delta \bu +(\lambda + \mu) \nabla(\div \bu)= \mathbf{f} \quad \quad  \bx \in \Omega
\]
where $\bu=\bu(\bx)$ is the displacement vector, $\Omega$ is the
spatial domain, and $\lambda$ and $\mu$ are the Lam\'e constants. A
mix of Dirichlet boundary conditions and traction conditions are applied to have a unique solution.
Discretization of the vector equation
leads to systems of equations whose coefficient matrix is s.p.d. and, since
each scalar component of the displacement vector is considered separately, has
a block form where each diagonal block corresponds to the matrix
coming from the discretization of Laplace equation for each unknown
component.
We considered Lam\'e equations on a beam characterized by $\mu=0.42$
and $\lambda=1.7$. The problem, which we refer to as $LE2D$, was
discretized using linear finite elements on triangular meshes of
three different sizes ($66690$, $264450$, $1053186$), obtained by uniform refinement using the software package
MFEM~\cite{mfem}.
\item[Parflow] A third set of s.p.d. linear systems comes from a groundwater model,
aimed at the numerical simulation of the filtration of 3D
incompressible single-phase flows through anisotropic porous media.
The linear systems arise from the discretization of the Darcy's equation,
with no-flow boundary conditions, performed by a cell-centered finite volume scheme
on a (structured) Cartesian grid. They were generated by using a Matlab
code implementing the fundamentals of reservoir simulations~\cite{A2007}
and can be regarded as simplified samples of systems arising in
ParFlow, a parallel computational model developed at the J\"{u}lich Supercomputing Centre (JSC).
We considered three different systems with dimension~$10^6$, corresponding to anisotropic
permeability tensors randomly generated
from a lognormal distribution having mean~1 and 5 different values of standard deviation, i.e., from 1 till 5, corresponding to increasing anisotropy levels, which we refer to as $ParflowN$, for $N=1, \ldots, 5$, respectively.
\end{description}

In all cases we solved the linear systems with right-hand sides set
equal to the unit vector.
The runs have been carried out on an Nvidia Titan V
(Nvidia Volta with 12 GB and 5120 CUDA cores running CUDA 9.1), operated by IAC-CNR in Rome.
Comparisons with single-node version of Nvidia AmgX 2.0.0.130-open source library have been carried out.

We always used the AMG preconditioner coupled with our flexible version of the
preconditioned CG solver (see Section~\ref{gpu}).
The CG procedure stopped when the euclidean norm of the relative
residual reached the tolerance $rtol=10^{-6}$ or the
number of iterations reached a predefined threshold $itmax=5000$.
We always considered AMG hierarchies,
with the maximum size of the coarsest matrix fixed to
$maxcoarseset *n^{1/3}$, where $n$ is the matrix dimension and
generally $maxcoarseset=40$.
A maximum number of levels was also fixed to $40$.
To maintain a maximum coarsening ratio equal to $4$,
we composed couples of prolongator operators computed by matching-based
aggregation, resulting in double pairwise aggregates. In all cases one sweep of $\ell^1-$Jacobi
method was applied as both pre- and post-smoother whereas 20 sweeps of the same method
were applied at the coarsest level.

\subsection{Performance Results}
\label{res1}

Hereafter we compare results obtained by applying the AMG preconditioner built by BootCMatchG  both as V-cycle and W-cycle.
Table~\ref{tabANI} summarizes performance results of the preconditioned CG, when V-cycle and W-cycle are applied, for the test cases $ANI1$ and $ANI2$ while system size increases. We report the execution time in milliseconds (ms), needed for the setup ($tsetup$) of the preconditioner, the number of iterations ($it$) and the time for the solution of the systems ($tsolve$) by the preconditioned flexible CG.
\begin{table}[t]
\centering{
{\begin{tabular}{||c||c|c||c|c||}
 \hline \hline
& \multicolumn{2}{c||}{\emph{V-cycle}} & \multicolumn{2}{c||}{\emph{W-cycle}} \\
 \hline \hline
 $tbuild $ & $it$ & $tsolve$ & $it$ & $tsolve$\\
 \hline  \hline
 \multicolumn{5}{||c||}{\emph{ANI1}} \\
 \hline \hline
 34.15  & 192 & 115.65  & 123 & 216.30\\
 64.60  & 302 & 464.69  & 161 & 664.84\\
 164.80 & 466 & 2476.52 & 205 & 2410.39\\
 \hline  \hline
 \multicolumn{5}{||c||}{\emph{ANI2}} \\
 \hline \hline
 30.34  & 194 & 116.66  & 123 & 213.90\\
 64.50  & 307 & 471.39  & 163 & 650.14\\
 167.59 & 481 & 2556.50 & 209 & 2437.63\\
 \hline \hline
\end{tabular}}
\caption{BootCMatchG: ANI test cases\label{tabANI}}
}
\end{table}
We observe that, as expected, W-cycle requires a smaller number of iterations than V-cycle, also showing a better algorithmic scalability. Indeed number of iterations increases more slowly for increasing problem size, for both test cases. On the other hand, W-cycle generally has a larger computational cost per iteration, therefore the best execution times are generally obtained by applying V-cycle, but, for the largest size, the large reduction of the number of iterations results also in a lower solving time for W-cycle.
Comparison between V-cycle and W-cycle are also reported, for the Parflow test cases, in Table~\ref{tabParflow}, for increasing level of anisotropy.
\begin{table}[t]
\centering{
{\begin{tabular}{||c||c|c||c|c||}
 \hline \hline
& \multicolumn{2}{c||}{\emph{V-cycle}} & \multicolumn{2}{c||}{\emph{W-cycle}} \\
 \hline \hline
 $tbuild $ & $it$ & $tsolve$ & $it$ & $tsolve$\\
 \hline  \hline
 \multicolumn{5}{||c||}{\emph{Parflow}} \\
 \hline \hline
 94.46  & 91 &  206.80 & 26  & 227.54\\
 96.06  & 82 & 189.98 & 42  & 390.43\\
 96.70 & 125 & 288.63 & 107  & 995.52\\
 96.10 & 87  & 198.36 & 27   & 236.69 \\
 97.83 & 604 & 1388.30 & 584 & 5569.89 \\
\hline \hline
\end{tabular}}
\caption{BootCMatchG: Parflow test cases\label{tabParflow}}
}
\end{table}
We observe that also for the \emph{Parflow} test cases, W-cycle requires a smaller number of iterations than V-cycle. On the other hand, in all cases V-cycle outperforms W-cycle in terms of execution times. Similar behaviour are observed also for the FE test cases, therefore,
we conclude that for the available choice of parallel smoother and coarsest solver, best execution times of BootCMatchG are generally obtained when V-cycle is applied.
\subsection{Comparison with AmgX}
\label{bcmvsamgx}

In the following we show a performance comparison with preconditioners implemented in the Nvidia AmgX package, when they run in a single-node setting.
We considered the two different configurations available for preconditioner setup in AmgX: classical AMG and aggregation-based AMG.
For classical AMG we used default configurations including \emph{D1}-interpolation and \emph{AHAT} strength of connection metric (see AmgX Reference Manual for details~\cite{amgxsw}), default parameters are also used for plain aggregation AMG, where aggregates of size $4$ are required. We refer to them as \emph{AmgXclassic} and \emph{AmgXaggr}, respectively, while the preconditioner implemented in BootCMatchG is referred to as \emph{BCMG}.
Both the preconditioner types are applied as V-cycle within preconditioned CG iterations and the same choices for pre/post-smoother and coarsest solver applied for BootCMatchG are considered for our comparisons. The parameters of Section~\ref{results} are also set for monitoring convergence of the preconditioned CG.

Note that we also applied the AmgX preconditioners as W-cycle. Furthermore, all the other available choices for smoothers and coarsest solvers have been considered for our test cases. In all cases, the best results in terms of execution times have been obtained with the AmgX preconditioners discussed in the present Section.

\begin{figure*}[htb]
\begin{center}
\includegraphics[width=0.45\textwidth]{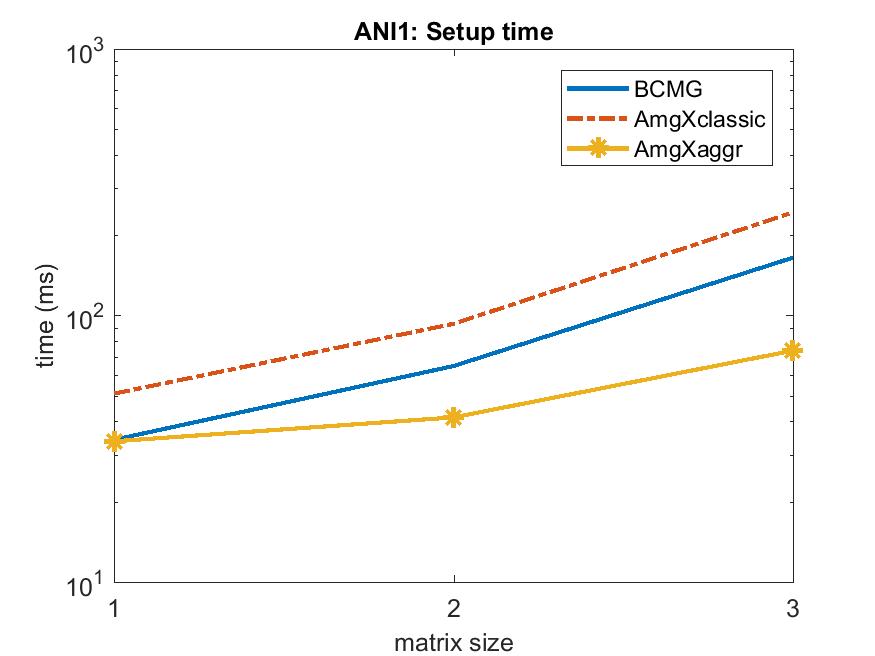}
\includegraphics[width=0.45\textwidth]{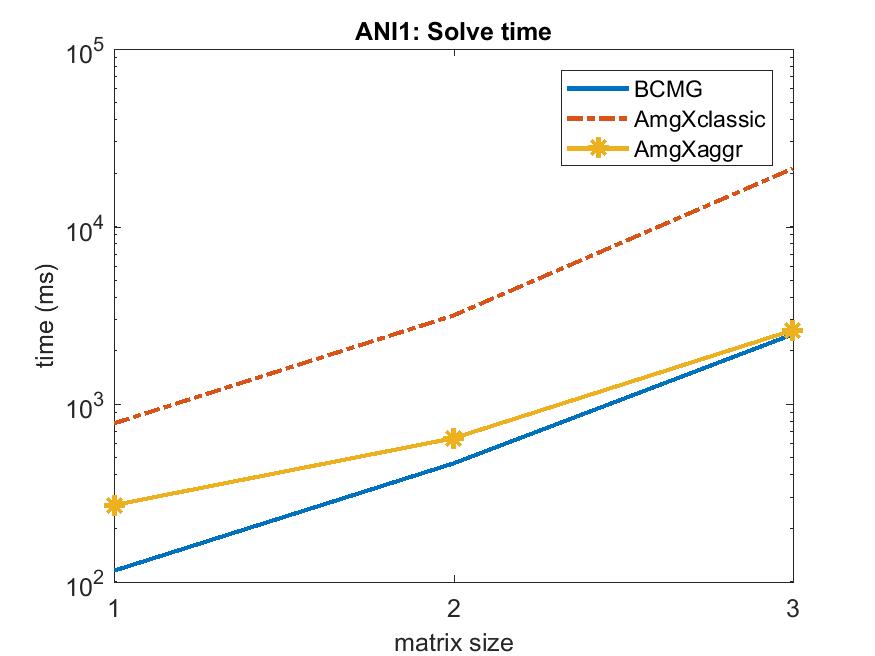}
\end{center}
\caption{ANI1 test case: AmgX vs BootCMatchG\label{fig-ANI1-1}}
\end{figure*}

\begin{figure*}[htb]
\begin{center}
\includegraphics[width=0.45\textwidth]{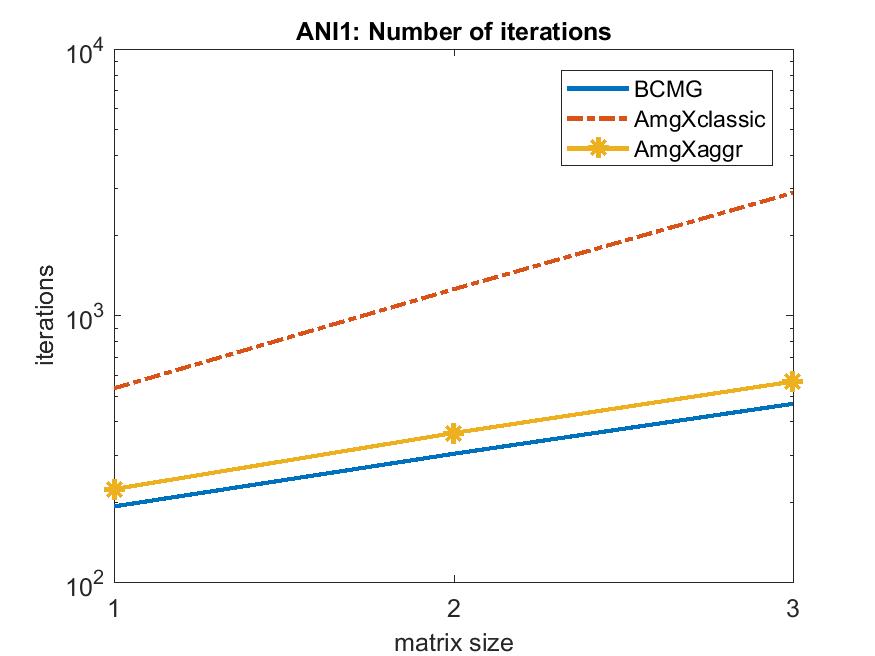}
\includegraphics[width=0.45\textwidth]{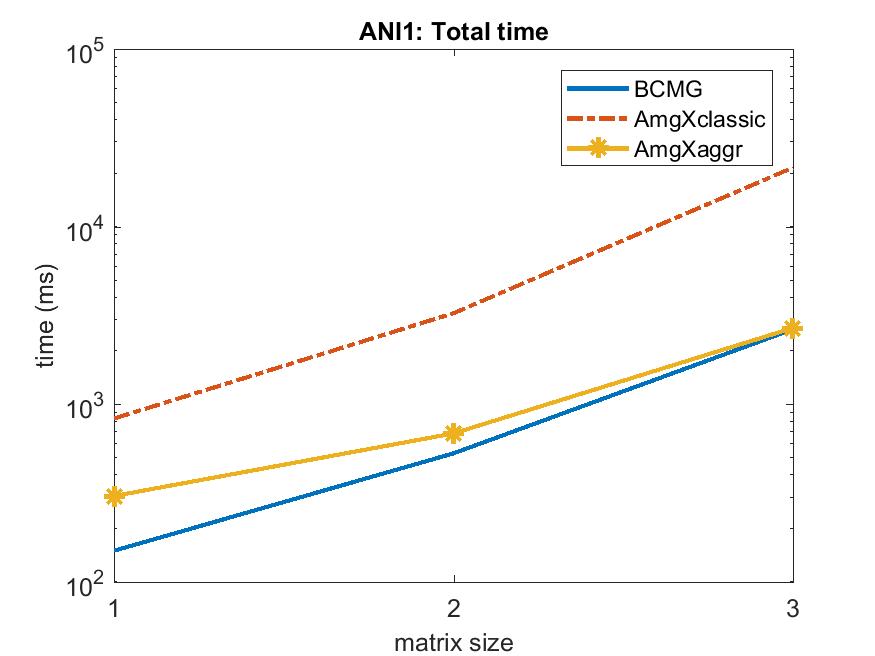}
\end{center}
\caption{ANI1 test case: AmgX vs BootCMatchG\label{fig-ANI1-2}}
\end{figure*}

In Fig.~\ref{fig-ANI1-1} we compare execution times for setup (left) and solve (right) phases of the preconditioners for the test case ANI1.
We observe that \emph{AmgXclassic} shows the longest times both for setup of the preconditioner and for solving the systems; its convergence behaviour
degrades significantly when the matrix size increases: the number of iterations increases from $533$ up to $2885$ going from the smallest to the largest size, as shown in Fig.~\ref{fig-ANI1-2} (left). Better numerical scalability properties are shown both for \emph{AmgXaggr} and for \emph{BCMG}, where a much more limited increase in the number of iterations is observed for increasing size. In particular, \emph{BCMG} generally shows the best behaviour in terms of number of iterations, which results in the best execution times for the solve phase and in total times better than or comparable with \emph{AmgXaggr}, as shown in Fig.~\ref{fig-ANI1-2} (right).

To gain a better understanding of our performance results, we analyzed some efficiency parameters of the various preconditioners. In Table~\ref{tabsetup-ANI1}, for increasing matrix size,
we summarize the number of levels $nl$ of the AMG preconditioner and the V-cycle operator complexity $Vcmplx=\frac{\sum_{k=1}^{nl}nnz(A^{k})}{nnz(A^1)}$, where $A^k$ is the matrix at level $k$ and $nnz(A^k)$ is the number of nonzeros of $A^k$, which give an estimate of the cost, in terms of both memory and computation requirements, of the preconditioner. We also report the average coarsening ratio of the AMG preconditioner
$cratio=\frac{1}{nl} \sum_{k=2}^{nl} \frac{n(A^{k-1})}{n(A^{k})}$, where $n(A^k)$ is the size of matrix $A^k$, which measures the ability of the coarsening schemes to obtain efficient AMG preconditioners with few levels and a limited operator complexity.
\begin{table*}[t]
\centering{
{\begin{tabular}{||c|c|c||c|c|c||c|c|c||}
 \hline \hline
\multicolumn{3}{||c||}{\emph{BCMG}} & \multicolumn{3}{|c||}{\emph{AmgXclassic}} & \multicolumn{3}{c||}{\emph{AmgXaggr}} \\
 \hline \hline
 $nl $ & $Vcmplx$ & $cratio$ & $nl $ & $Vcmplx$ & $cratio$ & $nl $ & $Vcmplx$ & $cratio$\\
 \hline  \hline
 \multicolumn{9}{||c||}{\emph{ANI1}} \\
 \hline \hline
 4 &  1.40  &  3.43 & 4 & 1.88 &  3.04 & 3 & 1.29 & 4.48\\
 5 &  1.40  &  3.11 & 5 & 1.89 &  3.02 & 4 & 1.30 & 4.48\\
 6 &  1.40  &  3.14 & 6 & 1.95 &  2.88 & 5 & 1.30 & 4.48\\
\hline \hline
\end{tabular}}
\caption{BootCMatchG vs AmgX: ANI1 test cases\label{tabsetup-ANI1}}
}
\end{table*}

We can observe that \emph{AMGXaggr} is able to build an AMG preconditioner with fewer levels and smaller operator complexities, with respect to \emph{BCMG} and \emph{AmgXclassic}, due to its ability to obtain larger coarsening ratios. This behaviour is the main reason of the shorter setup times of \emph{AMGXaggr}, indeed
  it requires one less coarsening step than the other preconditioners for all matrix sizes. On the other hand,
the quality of \emph{BCMG} appears better, indeed it requires fewer iterations for the PCG convergence leading to shorter solving times.
Similar results are obtained for the ANI2 test case, therefore, we omit them for sake of space.

\begin{figure*}[htb]
\begin{center}
\includegraphics[width=0.45\textwidth]{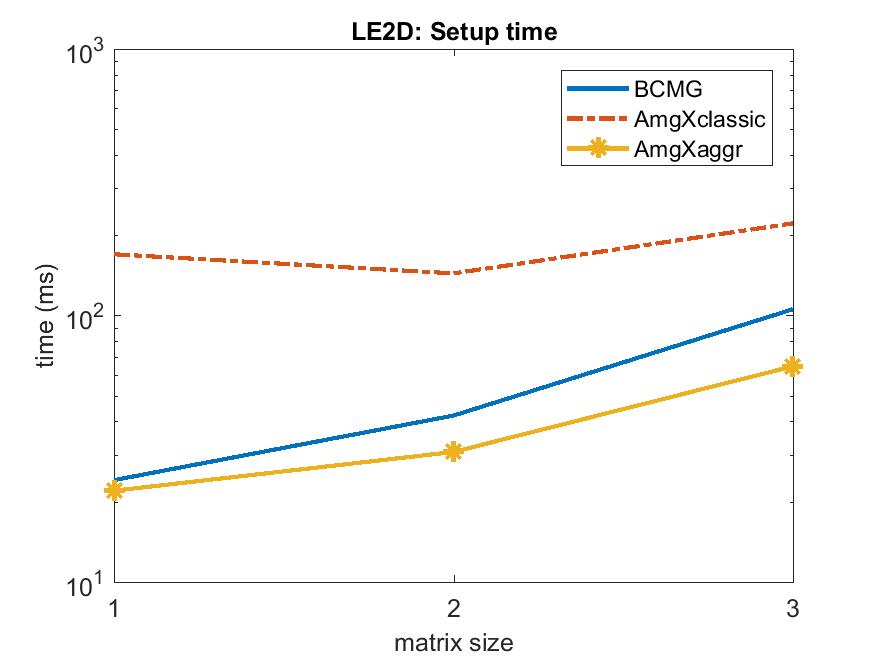}
\includegraphics[width=0.45\textwidth]{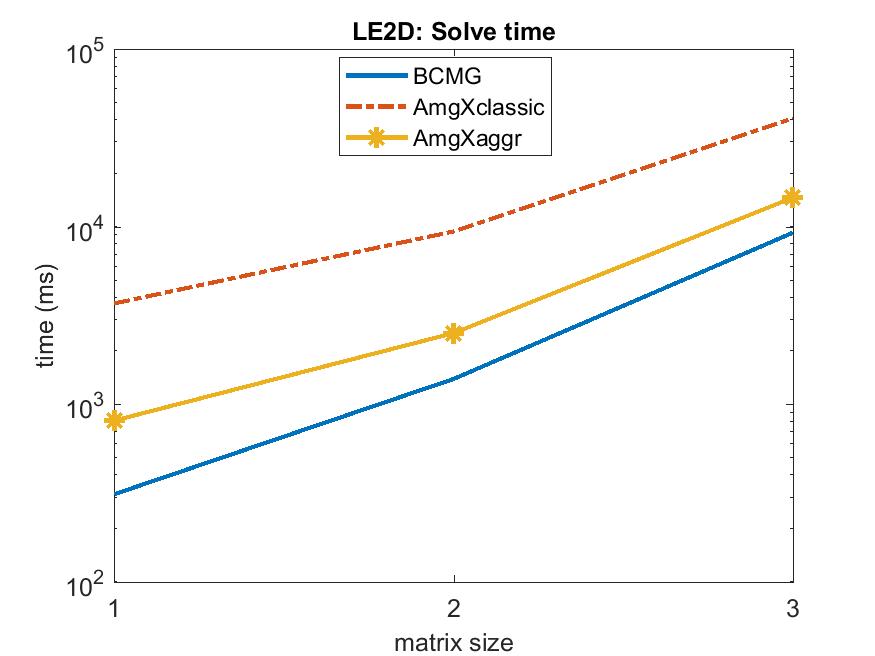}
\end{center}
\caption{LE2D test case: AmgX vs BootCMatchG\label{fig-LE2d-1}}
\end{figure*}

\begin{figure*}[htb]
\begin{center}
\includegraphics[width=0.45\textwidth]{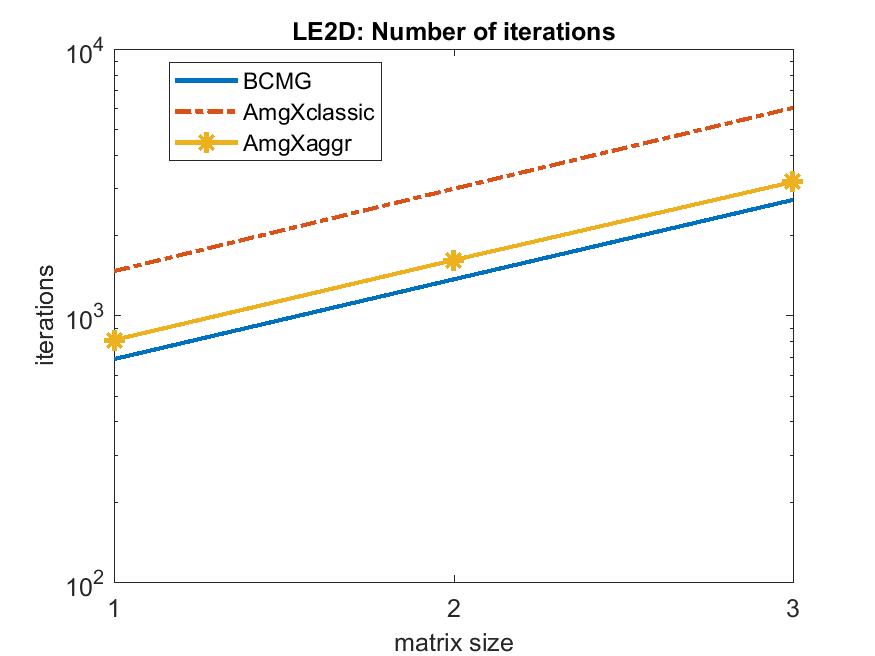}
\includegraphics[width=0.45\textwidth]{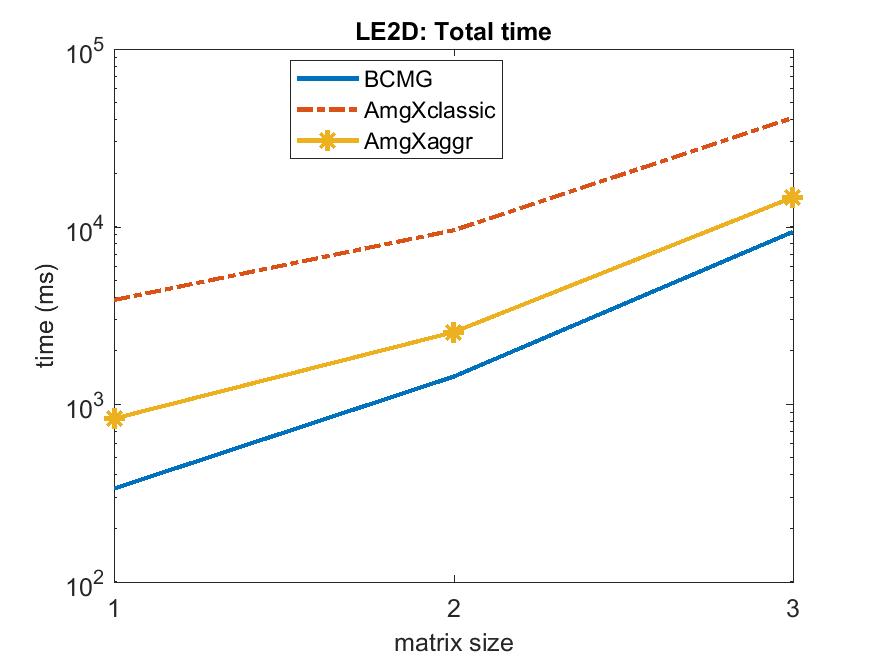}
\end{center}
\caption{LE2D test case: AmgX vs BootCMatchG\label{fig-LE2d-2}}
\end{figure*}

In Figs.~\ref{fig-LE2d-1}-\ref{fig-LE2d-2} we compare results obtained by the different preconditioners on the LE2D test case.
General behaviour is very similar to that obtained in the previous test case. We observe that, also in this case, \emph{AmgXaggr} shows
better scalability for the setup of the preconditioner, whereas \emph{BCMG} outperforms both AmgX preconditioners in the solve phase, due to the better convergence behaviour. Indeed, for all matrix sizes, \emph{BCMG} requires the smallest number of iterations. In this case, the significative reduction in the number of iterations is able to balance the longer setup times of \emph{BCMG}, resulting in the best total execution times.

In Table~\ref{tabsetup-LE2d}, we report efficiency parameters of the preconditioners for the LE test case. We note that in this case, \emph{BCMG} shows better averaged coarsening ratio than the ANI1 test case and it is able to obtain preconditioners with the same number of levels as \emph{AmgXaggr}. However, operator complexity of \emph{BCMG} is yet slightly larger, showing that coarse matrices are slightly more dense than \emph{AmgXaggr}. This behavior explains the slightly larger
setup times for \emph{BCMG}. On the other hand, the quality of our preconditioner is significative better, leading to good scalability and the best total execution times.
\begin{table*}[t]
\centering{
{\begin{tabular}{||c|c|c||c|c|c||c|c|c||}
 \hline \hline
\multicolumn{3}{||c||}{\emph{BCMG}} & \multicolumn{3}{|c||}{\emph{AmgXclassic}} & \multicolumn{3}{c||}{\emph{AmgXaggr}} \\
 \hline \hline
 $nl $ & $Vcmplx$ & $cratio$ & $nl $ & $Vcmplx$ & $cratio$ & $nl $ & $Vcmplx$ & $cratio$\\
 \hline  \hline
 \multicolumn{9}{||c||}{\emph{LE2D}} \\
 \hline \hline
 3 &  1.42  &  3.55 & 2 & 1.86 &  3.38 & 3 & 1.37 & 4.05\\
 4 &  1.44  &  3.51 & 3 & 1.90 &  3.63 & 4 & 1.38 & 3.96\\
 5 &  1.44  &  3.21 & 4 & 1.88 &  3.72 & 5 & 1.38 & 3.85\\
\hline \hline
\end{tabular}}
\caption{BootCMatchG vs AmgX: LE2D test cases\label{tabsetup-LE2d}}
}
\end{table*}

\begin{figure*}[htb]
\begin{center}
\includegraphics[width=0.45\textwidth]{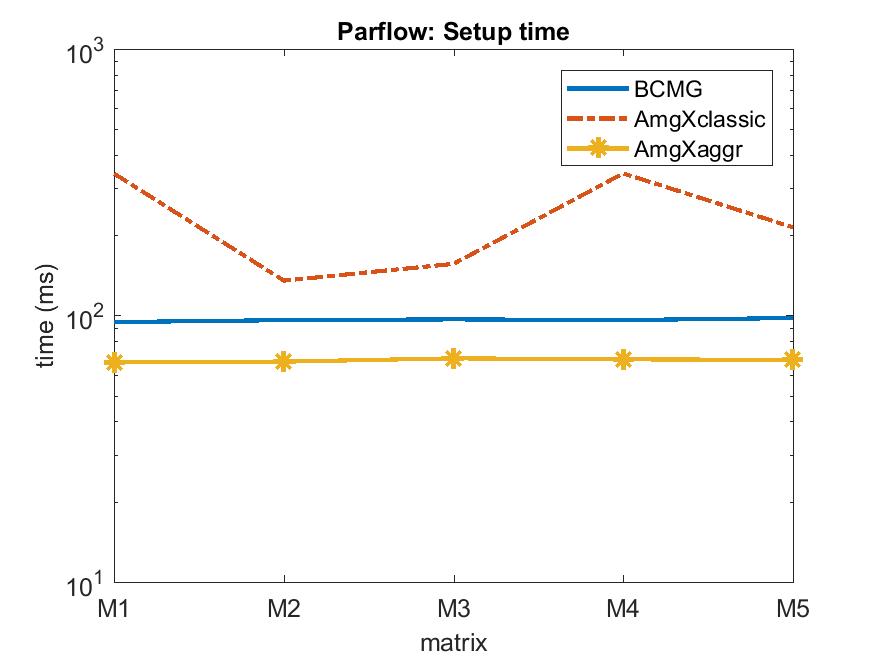}
\includegraphics[width=0.45\textwidth]{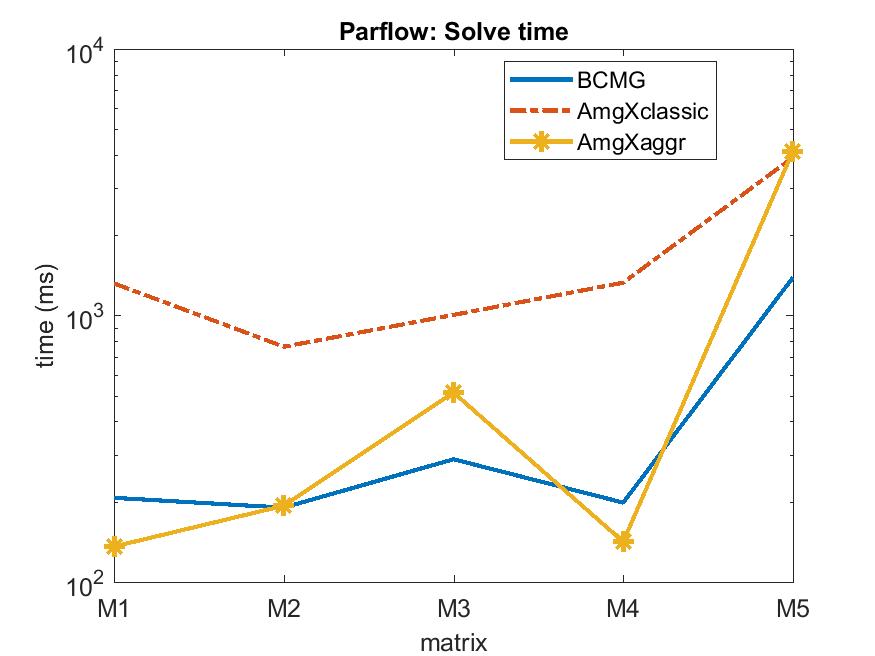}
\end{center}
\caption{Parflow test case: AmgX vs BootCMatchG\label{fig-Parflow-1}}
\end{figure*}

\begin{figure*}[htb]
\begin{center}
\includegraphics[width=0.45\textwidth]{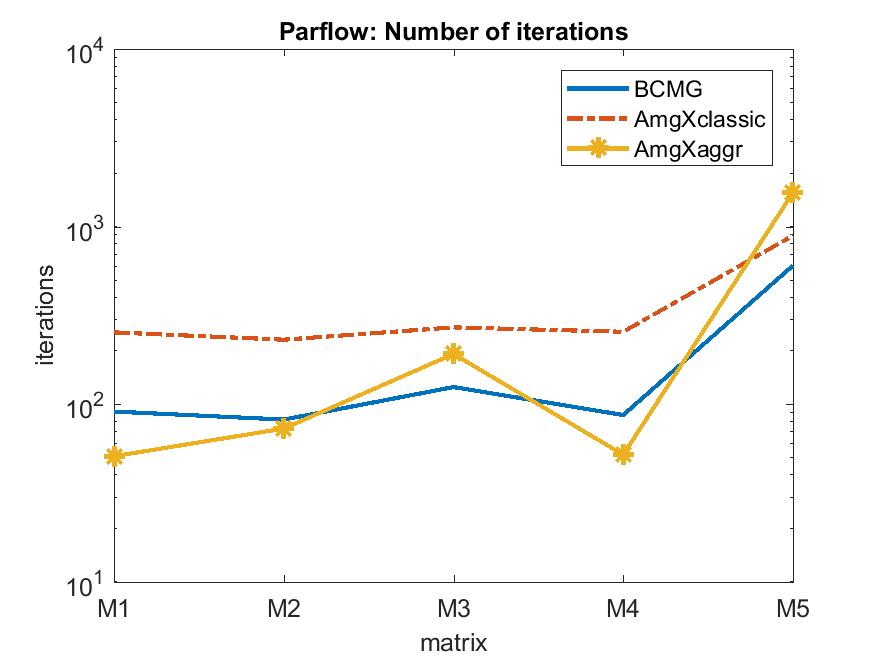}
\includegraphics[width=0.45\textwidth]{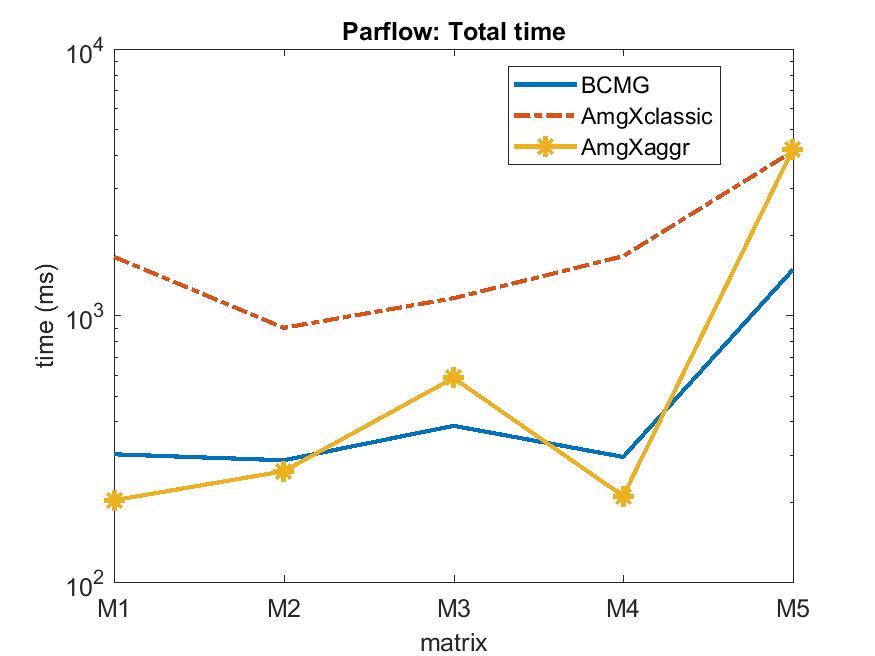}
\end{center}
\caption{Parflow test case: AmgX vs BootCMatchG\label{fig-Parflow-2}}
\end{figure*}

Finally, in Figs.~\ref{fig-Parflow-1}-\ref{fig-Parflow-2}, we report performance results of the preconditioners on the Parflow test cases.
We observe that, also in this case, \emph{AmgXclassic} has the worst behavior, both in the setup and in the solve phase.
In all cases \emph{AmgXaggr} has the best setup times. However, we note that while \emph{AmgXaggr} builds preconditioners with $4$ levels for all the Parflow test cases, \emph{BCMG}, due to smaller coarsening ratios, builds preconditioners with $6$ levels. This results in a better ratio
between setup times and number of levels for \emph{BCMG}, showing a good efficiency in the implementation of the basic parallel kernels of the setup phase of our preconditioner.
For Parflow1 and Parflow4, \emph{AmgXaggr} requires slightly fewer iterations than \emph{BCMG}, whereas \emph{BCMG} largely outperforms \emph{AmgXaggr} in the case of Parflow5, showing a better robustness with respect to anisotropy levels. We finally observe that \emph{BCMG} generally shows the best total execution time per iteration which ranges from $2.46$ ms of Parflow5 to $3.49$ ms of Parflow2.
\section{Conclusions}
\label{concl}

We presented the \emph{BootCMatchG} package, for
preconditioning and solving sparse s.p.d. linear systems on modern GPU
architectures.  The code implements an iterative linear solver of
Krylov type coupled with an AMG preconditioner based on the so-called
compatible weighted matching aggregation algorithm. We exploited
fine-grained parallelism and optimized global memory access in each
kernel both for the setup and the application of the AMG
preconditioner, as well as in the implementation of the Krylov solver.
We have rethought all main algorithms of the original package
available for standard CPU, by selecting and optimizing numerical
kernels for effective use of modern GPUs. To this aims, highly
parallel approximate matching algorithm and a robust version of the
Jacobi relaxation method were employed. Furthermore, we introduced the
concept of {\em miniwarp} for accessing GPU global memory and using
the available computing resources in an effective way.  We discussed
results for a large set of s.p.d. scalar and vector linear systems and
demonstrated that our solver outperforms the single-node Nvidia AmgX
library. Future work includes the exploitation of further parallel
smoothers, such as sparse approximate inverses, and a multi-GPU
version of the code.\\
The current version of the source code is available on request, by sending an email to one of the authors. In the near future we will make it available in a public repository ({\em e.g.,} github).

\section{Acknowledgements}
\label{ack}

The authors wish to thank Stefan Kollet and Wendy Sharples (JSC, J\"{u}lich) for making available the Matlab miniapp needed to generate Parflow test cases.

\bibliography{bibamggpu-finrel}

\end{document}